\documentclass{IEEEtran4PSCC}
\ifCLASSINFOpdf
   \usepackage[pdftex]{graphicx}
\else
   \usepackage[dvips]{graphicx}
\fi
\usepackage[cmex10]{amsmath}
\makeatletter
\let\old@ps@headings\ps@headings
\let\old@ps@IEEEtitlepagestyle\ps@IEEEtitlepagestyle
\def\psccfooter#1{%
    \def\ps@headings{%
        \old@ps@headings%
        \def\@oddfoot{\strut\hfill#1\hfill\strut}%
        \def\@evenfoot{\strut\hfill#1\hfill\strut}%
    }%
    \def\ps@IEEEtitlepagestyle{%
        \old@ps@IEEEtitlepagestyle%
        \def\@oddfoot{\strut\hfill#1\hfill\strut}%
        \def\@evenfoot{\strut\hfill#1\hfill\strut}%
    }%
    \ps@headings%
}
\makeatother

\usepackage{tikz,pgfplots}
\usepackage{algorithm}
\usepackage{algpseudocode}
\usepackage{booktabs}
\usepackage{amssymb}
\usepackage{enumitem}

\begin{document}
%
\title{The Sweet Spot of Bound Tightening \\ for Topology Optimization}

\author{
\IEEEauthorblockN{Salvador Pineda, Juan Miguel Morales}
\IEEEauthorblockA{OASYS research group, University of Málaga, Spain \\
\{spineda, juan.morales\}@uma.es}
}


\maketitle

\begin{abstract}
Topology optimization has emerged as a powerful and increasingly relevant strategy for enhancing the flexibility and efficiency of power system operations. However, solving these problems is computationally demanding due to their combinatorial nature and the use of big-M formulations. Optimization-based bound tightening (OBBT) is a well-known strategy to improve the solution of mixed-integer linear programs (MILPs) by computing tighter bounds for continuous variables. Yet, existing OBBT approaches in topology optimization typically relax all switching decisions in the bounding subproblems, leading to excessively loose feasible regions and limited bound improvements. In this work, we propose a topology-aware bound tightening method that uses network structure to determine which switching variables to relax. Through extensive computational experiments on the IEEE 118-bus system, we find that keeping a small subset of switching variables as binary, while relaxing the rest, strikes a sweet spot between the computational effort required to solve the bounding problems and the tightness of the resulting bounds.
\end{abstract}


\thanksto{\noindent The work of S. Pineda and J. M. Morales was supported by the Spanish Ministry of Science and Innovation (AEI/10.13039/501100011033) through project  PID2023-148291NB-I00 and by the Consejería de Universidad, Investigación e Innovación de la Junta de Andalucía through FEDER funds (grant PPRO-TEP967-G-2023). S. Pineda, and J. M. Morales are with the research group OASYS, University of Malaga, Malaga 29071, Spain:  spineda@uma.es; juan.morales@uma.es. Finally, the authors thankfully acknowledge the computer resources, technical expertise, and assistance provided by the SCBI (Supercomputing and Bioinformatics) center of the University of M\'alaga.}

\section{Introduction} \label{sec:intro}

Topology optimization has emerged as a powerful and increasingly relevant strategy for enhancing the flexibility and efficiency of power system operations. This flexibility is particularly valuable in systems with high shares of renewable generation, where variability and uncertainty pose significant operational challenges \cite{o2005dispatchable}. Topology switching can help relieve congestion, reduce operational costs, and improve the accommodation of renewables, all while maintaining reliability. When this flexibility is used in a systematic and optimized manner, the problem is referred to as Optimal Transmission Switching (OTS) \cite{fisher2008optimal}. In the general case, OTS seeks to determine the best combination of line statuses (on/off) to minimize operational costs subject to power flow and reliability constraints. 

The OTS problem is inherently combinatorial due to binary line status decisions and remains challenging even under DC approximations. The DC-OTS problem is typically formulated as a MILP using big-M constants \cite{pineda2024tight}. While such formulations are solver-friendly, their performance is highly sensitive to the choice of big-M values: if too large, they yield weak relaxations and high computational cost; if too small, they may exclude feasible solutions. Hence, computing valid and tight big-M values is crucial for solving OTS efficiently. 

Although several methods have been proposed to determine these constants (see \cite{numan2023role}), many lack general validity as demonstrated in \cite{fattahi2019bound} using counterexamples. To address this, the authors of \cite{fattahi2019bound} propose a method to compute valid big-M values when the transmission network includes a \emph{spanning tree of non-switchable lines}, by solving shortest-path problems. These values are easy to obtain and valid for any operating condition, though typically not very tight. To improve tightness, \cite{pineda2024tight} propose an optimality-based bound tightening (OBBT) approach that solves a pair of optimization problems per line, assuming a non-switchable spanning tree. Although more computationally demanding, this method produces tighter big-M values adapted to the current system operating conditions and significantly accelerates the overall OTS solution.

For the more challenging case where \emph{all lines are switchable}, Moulin et al.~\cite{moulin2010transmission} show that valid big-M values can, in principle, be computed based on the longest path between node pairs in the network. However, computing such paths is itself a combinatorially hard problem and they propose a conservative and easily computable upper bound as a surrogate. While this approximation guarantees validity, it is often excessively conservative and leads to weak relaxations that limit computational efficiency. One could also apply the OBBT technique proposed in \cite{pineda2024tight}, but its effectiveness is severely limited in the fully switchable scenario. In this case, the bounding problems solved within OBBT lead to overly large and weak feasible regions, resulting in bounds that are often too loose to provide any meaningful tightening. As a result, the main difficulty in this setting stems not only from the increased number of binary variables, but more critically, from the challenge of computing sufficiently tight and valid big-M values.

In this paper, we tackle the challenge of solving the Optimal Transmission Switching (OTS) problem in fully switchable networks by proposing a topology-aware bound tightening method to enhance the big-M constants used in its MILP formulation. Rather than relaxing all switching decisions, our approach retains a selected subset of binary variables (chosen based on network topology) while relaxing the rest. This yields bounding problems with a more meaningful feasible region and leads to tighter and more informative big-M values, even when all lines are switchable. Our main contributions are:
\begin{itemize}
\item We introduce a tunable, topology-aware method that retains integrality for a subset of switching variables during bound tightening. This partial relaxation preserves relevant discrete structure, improving bound accuracy while maintaining computational tractability.
\item We demonstrate, through extensive experiments, that this strategy strikes an effective balance between bound quality and complexity, leading to faster and more reliable solutions to the overall OTS problem.
\end{itemize}

The remainder of the paper is organized as follows. Section~\ref{sec:ots} introduces the Optimal Transmission Switching problem based on linear DC power flow equations, presents its reformulation as a mixed-integer linear program, and discusses the impact of big-M bounds on its computational complexity. In Section~\ref{sec:boundtightening}, we first describe the classical optimization-based bound tightening technique and then present our topology-aware proposal. Section~\ref{sec:comparisson} introduces the benchmark approaches used for comparison, along with the performance metrics employed to assess solution quality and computational effort. Computational results for the IEEE 118-bus system are discussed in Section~\ref{sec:results}. Finally, conclusions are duly drawn in Section~\ref{sec:conclusions}.

\section{Optimal Transmission Switching} \label{sec:ots}

The Optimal Transmission Switching (OTS) problem can be formulated either using a full \emph{Alternating Current} (AC) power flow model or using the simpler \emph{Direct Current} (DC) approximation. For readers interested in a broader discussion of AC versus DC formulations and their respective challenges, we refer to review works such as \cite{numan2023role,  flores2020alternative, crozier2022feasible}. While the AC-OTS formulation captures voltage magnitudes, reactive power, and line losses, and ensures solutions that are AC-feasible, it is considerably more complex to solve. Nevertheless, several studies have addressed AC-OTS directly \cite{soroush2013accuracies,  barrows2014correcting, kocuk2017new}, demonstrating its importance for ensuring system security and accurate cost assessment. At the same time, a large body of literature continues to adopt the DC-OTS formulation \cite{fisher2008optimal, kocuk2016cycle, fattahi2019bound, liu2012heuristic, barrows2012computationally, fuller2012fast, hinneck2022optimal, johnson2020knearest, pineda2024learning, yang2019line, dey2022node, hedman2012flexible, moulin2010transmission}, which provides a tractable and widely accepted benchmark for testing solution approaches. For this reason, and to focus on the analysis of bound-tightening and reformulation strategies, we consider the DC-OTS formulation in this work.

Consider a power system represented by a set of nodes $\mathcal{N}$ and transmission lines $\mathcal{L}$. For simplicity, we assume that each node $n \in \mathcal{N}$ is equipped with a generator and a load. Let $p_n$ and $d_n$ denote the generation and demand at node $n$, respectively. Each generator operates within its technical limits, defined by the minimum and maximum power outputs $\underline{p}_n$ and $\overline{p}_n$, and incurs a marginal cost $c_n$ per unit of energy produced. Each transmission line $l = (n,m) \in \mathcal{L}$ connects nodes $n$ and $m$, and supports a power flow denoted by $f_l$. By convention, $f_l > 0$ indicates flow from node $n$ to node $m$, and $f_l < 0$ indicates flow in the opposite direction. The power flow on line $l$ is constrained by its thermal limits, with $\underline{f}_l<0$ and $\overline{f}_l>0$ denoting its minimum and maximum allowable flows, respectively. For modeling purposes, we define a dummy power flow variable $\tilde{f}_l$ for each line $l \in \mathcal{L}$, given by the product of the line's susceptance $b_l$ and the voltage angle difference across its terminal nodes: $\tilde{f}_l = b_l(\theta_n - \theta_m)$. Since all lines are considered switchable, their operational status is represented by a binary variable $x_l$. If $x_l = 1$, the line is active and the actual power flow is $f_l = \tilde{f}_l$; if $x_l = 0$, the line is disconnected and the power flow is zero, i.e., $f_l = 0$. With this notation in place, the DC-OTS problem can be formulated as follows:
\begin{subequations}\label{eq:OTS_NP}
\begin{IEEEeqnarray}{l}
\min \quad \sum_{n} c_{n} \, p_{n} \label{eq:OTS_NP_obj}\\
\text{subject to}  \nonumber \\ 
f_l = x_l \tilde{f}_l, \quad \forall l \in \mathcal{L} \label{eq:OTS_NP_Flow}\\
\underline{f}_l \leq f_l \leq \overline{f}_l, \quad \forall l \in \mathcal{L} \label{eq:OTS_NP_Flow_limit_S}\\
\tilde{f}_l = b_l(\theta_n-\theta_m), \quad \forall l=(n,m) \in \mathcal{L} \label{eq:OTS_NP_Flow_Dummy}\\
p_n - d_n = \sum_{l\in\mathcal{L}(n,\cdot)} f_l - \sum_{l\in\mathcal{L}(\cdot,n)} f_l, \quad \forall n \in \mathcal{N} \label{eq:OTS_NP_PB}\\
\underline{p}_n \leq p_n \leq \overline{p}_n, \quad \forall n \in \mathcal{N} \label{eq:OTS_NP_Plimits}\\
\theta_1 = 0 \label{eq:OTS_NP_slack}\\
x_l \in \{0,1\}, \quad \forall l \in \mathcal{L} \label{eq:OTS_NP_binary}
\end{IEEEeqnarray}
\end{subequations}

The objective function in equation~\eqref{eq:OTS_NP_obj} minimizes the total generation cost across the network. Power flows along transmission lines are modeled by equation~\eqref{eq:OTS_NP_Flow} and its limits are imposed by constraints~\eqref{eq:OTS_NP_Flow_limit_S}. The dummy power flow variable is defined in equation~\eqref{eq:OTS_NP_Flow_Dummy}, the nodal power balance is enforced through equation~\eqref{eq:OTS_NP_PB} and the generator outputs are constrained by their operating limits, as specified in equation~\eqref{eq:OTS_NP_Plimits}. The voltage angle at a reference node is fixed to zero in equation~\eqref{eq:OTS_NP_slack}. Finally, the binary nature of the switching variables is enforced by constraint~\eqref{eq:OTS_NP_binary}. Problem~\eqref{eq:OTS_NP} is a mixed-integer nonlinear programming problem due to the product $x_{l}\tilde{f}_{l}$ in~\eqref{eq:OTS_NP_Flow}. However, constraint~\eqref{eq:OTS_NP_Flow} can be linearized by introducing a pair of large enough constants $\underline{M}_{l}<0$, $\overline{M}_{l}>0$ per line \cite{hedman2012flexible} as follows:
\begin{subequations}\label{eq:linear}
\begin{align}
& (1-x_l)\underline{M}_l \leq -f_l + \tilde{f}_l \leq (1-x_l)\overline{M}_{l} \\
& x_l\underline{f}_l \leq f_l \leq x_l\overline{f}_l
\end{align}
\end{subequations}
where $\underline{M}_{l}, \overline{M}_{l}$ are guaranteed to be valid bounds of the dummy flow variable $\tilde{f}_l$ when the line $l$ is disconnected ($x_l = 0$). Under that assumption, the DC-OTS is reformulated as the following mixed-integer linear programming problem
\begin{equation} \label{eq:ots_mip}
\min  \quad  \sum_{n} c_{n} \, p_{n} \quad
\text{subject to}  \quad \eqref{eq:OTS_NP_Flow_Dummy}-\eqref{eq:OTS_NP_binary}, \eqref{eq:linear}    
\end{equation} 

While model~\eqref{eq:ots_mip} can be solved with general-purpose solvers like Gurobi~\cite{gurobi}, its performance depends critically on the tightness of the big-M parameters $\underline{M}_l$ and $\overline{M}_l$ \cite{camm1990cutting, crema2014mathematical}. Loose values lead to weak relaxations and poor solver performance, as previously discussed. The next section describes how these parameters are typically initialized and how they can be systematically refined to enhance computational efficiency.

Finally, we note that the formulation above corresponds to a \emph{base-case} DC-OTS problem and does not explicitly enforce N--1 security constraints, which are typically required in operational settings. Incorporating N--1 security would require extending the model to a security-constrained OTS formulation, where power balance and line flow constraints must be satisfied not only in the nominal topology but also under a predefined set of single-component contingencies. While such extensions are well studied in the context of security-constrained OPF, their integration with transmission switching significantly increases both the size and the complexity of the resulting mixed-integer problem. Since the primary goal of this work is to analyze and improve the computational performance of bound-tightening techniques for DC-OTS formulations, we deliberately focus on the base-case problem to isolate the impact of the proposed methodology.

\section{Bound Tightening}
\label{sec:boundtightening}

In this section, we describe an optimization-based bound tightening procedure that improves upon the approach presented in~\cite{pineda2024tight} for the case where all lines are switchable. To concisely formulate the methodology proposed in~\cite{pineda2024tight}, we define the following vectors of decision variables: $\mathbf{p} = [p_n,\ n \in \mathcal{N}]$, $\boldsymbol{\theta} = [\theta_n,\ n \in \mathcal{N}]$, $\mathbf{f} = [f_l,\ l \in \mathcal{L}]$, $\tilde{\mathbf{f}} = [\tilde{f}_l,\ l \in \mathcal{L}]$, and $\mathbf{x} = [x_l,\ l \in \mathcal{L}]$. The corresponding parameter vectors are defined as $\mathbf{F} = [(\underline{f}_l, \overline{f}_l),\ l \in \mathcal{L}]$ and $\mathbf{M} = [(\underline{M}_l, \overline{M}_l),\ l \in \mathcal{L}]$. The initial power flow bounds are denoted as $\mathbf{F}^0$ and just correspond to the line thermal limits. The initial big-M values, also denoted as $\mathbf{M}^0$, are computed using the methodology proposed in \cite{moulin2010transmission}. For given values of $\mathbf{F}$ and $\mathbf{M}$, we denote by $\mathcal{R}(\mathbf{F}, \mathbf{M})$ the feasible region defined by equations~\eqref{eq:OTS_NP_Flow_Dummy}--\eqref{eq:OTS_NP_slack} and~\eqref{eq:linear}. To represent the space of binary and relaxed line statuses, we define $\mathcal{X}^B := \{\mathbf{x} \in \{0,1\}^{|\mathcal{L}|}\}$ and its relaxed counterpart $\mathcal{X}^R := \{\mathbf{x} \in [0,1]^{|\mathcal{L}|}\}$. Additionally, we define the set $\mathcal{C} := \{\mathbf{p} \in \mathbb{R}^{|\mathcal{N}|} : \sum_n c_n p_n \leq \overline{C} \}$, where $\overline{C}$ is an upper bound on the optimal generation cost of problem~\eqref{eq:ots_mip}. With this notation in place, the methodology in~\cite{pineda2024tight} updates the bounds $\underline{f}_l$, $\overline{f}_l$, $\underline{M}_l$, and $\overline{M}_l$ for a given line $l \in \mathcal{L}$ by solving the following four bounding problems:
\begin{subequations} \label{eq:bounding}
\begin{align}
& \underline{f}_l= \underset{\mathcal{R}(\mathbf{F},\mathbf{M})\,\cap\,\mathcal{X}^R\,\cap\,\mathcal{C}}{\min} \quad f_l \quad \text{s.t.} \quad x_l = 1 \label{eq:bounding_min_flow}\\
& \overline{f}_l= \underset{\mathcal{R}(\mathbf{F},\mathbf{M})\,\cap\,\mathcal{X}^R\,\cap\,\mathcal{C}}{\max} \quad f_l \quad \text{s.t.} \quad x_l = 1 \label{eq:bounding_max_flow}\\
& \underline{M}_l= \underset{\mathcal{R}(\mathbf{F},\mathbf{M})\,\cap\,\mathcal{X}^R\,\cap\,\mathcal{C}}{\min} \quad \tilde{f}_l \quad \text{s.t.} \quad x_l = 0 \label{eq:bounding_min_dummy}\\
& \overline{M}_l= \underset{\mathcal{R}(\mathbf{F},\mathbf{M})\,\cap\,\mathcal{X}^R\,\cap\,\mathcal{C}}{\max} \quad \tilde{f}_l \quad \text{s.t.} \quad x_l = 0 \label{eq:bounding_max_dummy}
\end{align}
\end{subequations}

The bounding problems in~\eqref{eq:bounding} are linear programs that can be solved efficiently. When the network contains a spanning tree of connected (non-switchable) lines, \cite{pineda2024tight} proves that these bounds can be significantly tightened with respect to those computed in \cite{fattahi2019bound}, enabling the resulting OTS problem to be solved more efficiently. In this paper, we focus on the more challenging case in which \emph{all} transmission lines are switchable. In this setting, the relaxed feasible region used in the bounding problems may become overly permissive, as all binary variables are allowed to take continuous values in $[0,1]$. As a result, the bounding problems in~\eqref{eq:bounding} may fail to produce any meaningful improvement over the initial bounds computed using the method proposed in~\cite{moulin2010transmission}.

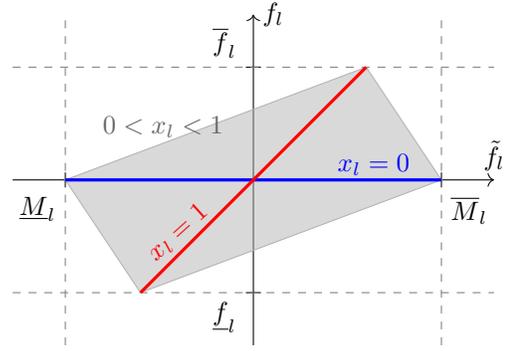
\begin{figure}
    \centering
    \begin{tikzpicture}
    \def\fMin{-1.5}    
    \def\fMax{1.5}     
    \def\ftildeMin{-2.5} 
    \def\ftildeMax{2.5}  

    \def\ftildeAxisMin{\ftildeMin - 0.7} 
    \def\ftildeAxisMax{\ftildeMax + 0.7} 

    \def\fAxisMin{\fMin - 0.7} 
    \def\fAxisMax{\fMax + 0.7} 

    \draw[->] (\ftildeAxisMin, 0) -- (\ftildeAxisMax, 0) node[above] {$\tilde{f}_l$}; 
    \draw[->] (0, \fAxisMin) -- (0, \fAxisMax) node[right] {$f_l$}; 

    \draw (\ftildeMin, 0.1) -- (\ftildeMin, -0.1) node[below left] {$\underline{M}_l$};
    \draw (\ftildeMax, 0.1) -- (\ftildeMax, -0.1) node[below right] {$\overline{M}_l$};

    \draw (0.1, \fMin) -- (-0.1, \fMin) node[below left] {$\underline{f}_l$};
    \draw (0.1, \fMax) -- (-0.1, \fMax) node[above left] {$\overline{f}_l$};

    \draw[dashed, gray] (\ftildeMin, \fAxisMin) -- (\ftildeMin, \fAxisMax);
    \draw[dashed, gray] (\ftildeMax, \fAxisMin) -- (\ftildeMax, \fAxisMax);

    \draw[dashed, gray] (\ftildeAxisMin, \fMin) -- (\ftildeAxisMax, \fMin);
    \draw[dashed, gray] (\ftildeAxisMin, \fMax) -- (\ftildeAxisMax, \fMax);    

    \pgfmathsetmacro{\polyARedPointAx}{\ftildeMin}
    \pgfmathsetmacro{\polyARedPointAy}{0}
    \pgfmathsetmacro{\polyARedPointBx}{\fMax}
    \pgfmathsetmacro{\polyARedPointBy}{\fMax}
    \pgfmathsetmacro{\polyARedPointCx}{\ftildeMax}
    \pgfmathsetmacro{\polyARedPointCy}{0}
    \pgfmathsetmacro{\polyARedPointDx}{\fMin}
    \pgfmathsetmacro{\polyARedPointDy}{\fMin}

    \filldraw[fill=black!30, draw=black!30, fill opacity=0.5]
        (\polyARedPointAx, \polyARedPointAy) -- (\polyARedPointBx, \polyARedPointBy) -- (\polyARedPointCx, \polyARedPointCy) -- (\polyARedPointDx, \polyARedPointDy) -- cycle;

    \draw[blue, very thick] (\ftildeMin, 0) -- (\ftildeMax, 0);   
    \draw[red, very thick] (\fMin, \fMin) -- (\fMax, \fMax);
    \node[red, rotate=42] at (-1, -0.7) {\textbf{$x_l=1$}};
    \node[blue] at (1.6, 0.2) {\textbf{$x_l=0$}};   
    \node[black!60] at (-1.2, 0.7) {$0 < x_l < 1 $};

\end{tikzpicture}    
    \caption{Relation between power flow variable $f_l$ and dummy variable $\tilde{f}_l$}
    \label{fig:relaxation}
\end{figure}

To illustrate the effect of relaxing all switching decisions in the bounding problems~\eqref{eq:bounding}, Fig.~\ref{fig:relaxation} plots the dummy variable $\tilde{f}_l$ on the x-axis and the actual power flow $f_l$ on the y-axis. Dashed lines indicate the upper and lower bounds for both variables. The bold red line corresponds to the case $x_l = 1$, where the line is fully operational: the power flow $f_l$ exactly matches the dummy variable $\tilde{f}_l$, which represents the product of the line susceptance and the angle difference, i.e., Kirchhoff’s law is enforced. The bold blue line shows the opposite case, $x_l = 0$, where the line is disconnected: $f_l = 0$, while $\tilde{f}_l$ remains unconstrained within its bounds. The gray area represents intermediate, fractional values of $x_l$ with $0 < x_l < 1$, and highlights two important effects. First, any fractional switching value leads to a reduction in the feasible flow region, i.e., $f_l$ lies strictly within its bounds, and the effective flow that can be transferred is reduced. Second, and more importantly, fractional values decouple $f_l$ and $\tilde{f}_l$, breaking the physical relationship implied by Kirchhoff’s law. As a result, the model behaves more like a generic flow network, where power can be rerouted without obeying physical constraints, akin to pipeline flow models.

Therefore, since the bounding problems proposed in~\cite{pineda2024tight} are solved by relaxing all switching decisions, the resulting power flows are not required to obey Kirchhoff’s law, which may lead to overly loose bounds. To address this issue, we propose solving the bounding problems over a less relaxed feasible region, in which a subset of the switching variables is kept binary while the rest are relaxed. For a given line $l \in \mathcal{L}$ and a closeness parameter $k \in \mathbb{N}$, we define the set $\mathcal{L}^k_l$ as the collection of lines that are within topological distance $k$ from line $l$, excluding $l$ itself. Specifically, the level-1 set $\mathcal{L}^1_l$ includes all lines in $\mathcal{L} \setminus \{l\}$ that are directly connected to either of the two end nodes of line $l$. The level-2 set $\mathcal{L}^2_l$ includes all lines in $\mathcal{L} \setminus \{l\}$ that are directly connected to any node incident to a line in $\mathcal{L}^1_l$. More generally, the sets are constructed recursively such that $\mathcal{L}^{k}_l$ contains all lines that are adjacent to nodes connected by any line in $\mathcal{L}^{k-1}_l$, satisfying the nesting property $\mathcal{L}^{1}_l \subseteq \mathcal{L}^{2}_l \subseteq \dots$. For completeness, we define $\mathcal{L}^0_l := \emptyset$.

Using the previously defined topology-based neighborhood, we now formalize the relaxation scheme applied to the binary variables. For a given line $l \in \mathcal{L}$ and a closeness level $k \in \mathbb{N}$, we define the set $\mathcal{X}^k_l$ as the space of line status vectors in which only the variables corresponding to lines in $\mathcal{L}^k_l$ are constrained to be binary, while the rest are relaxed. Formally,
\begin{equation}\label{eq:X_k_l}
\mathcal{X}^k_l := \left\{ \mathbf{x} \in [0,1]^{|\mathcal{L}|} :  x_j \in \{0,1\} \ \forall j \in \mathcal{L}^k_l \right\}
\end{equation}

We also formulate the new bounding problems~\eqref{eq:bounding_new} using the partially relaxed feasible set $\mathcal{X}^k_l$ as follows:
\begin{subequations} \label{eq:bounding_new}
\begin{align}
& \underline{f}_l\,/\,\overline{f}_l= \underset{\mathcal{R}(\mathbf{F},\mathbf{M})\,\cap\,\mathcal{X}^k_l\,\cap\,\mathcal{C}}{\min/\max} \quad f_l \quad \text{s.t.} \quad x_l = 1 \label{eq:bounding_new_min_flow}\\
& \underline{M}_l\,/\,\overline{M}_l= \underset{\mathcal{R}(\mathbf{F},\mathbf{M})\,\cap\,\mathcal{X}^k_l\,\cap\,\mathcal{C}}{\min/\max} \quad \tilde{f}_l \quad \text{s.t.} \quad x_l = 0 \label{eq:bounding_new_min_dummy}
\end{align}
\end{subequations}

The intuition behind the construction of the set~\eqref{eq:X_k_l} is that by enforcing integrality for lines in the neighborhood of $l$, we preserve the physical consistency of power flows in that area, i.e., we ensure that Kirchhoff's laws are satisfied locally. In contrast, switching decisions for lines located farther away are relaxed, allowing greater flexibility and reducing computational burden. This selective relaxation introduces a controllable trade-off between tractability and bound quality. As the parameter $k$ increases, the feasible region becomes progressively less relaxed, yielding tighter bounds but requiring more computational resources to solve.

The overall procedure to initialize the bounds, perform tightening, and solve the OTS formulation is summarized in Algorithm~\ref{alg:ots_bound_tightening}. We refer to the proposed method as \texttt{TBT-$k$}, standing for \emph{Topological Bound Tightening} with level-$k$ partial relaxation. For completeness, note that \texttt{TBT-0} corresponds exactly to the bound tightening approach originally proposed in~\cite{pineda2024tight}, in which all binary variables are relaxed during the bounding process.

\begin{algorithm}
\begin{small}
\caption{Bound Tightening and OTS Solution Procedure}
\label{alg:ots_bound_tightening}
\begin{algorithmic}
\State \textbf{Input:} Closeness level $k$.

\begin{enumerate}[leftmargin=15pt, label={\arabic*)}]
    \item Compute initial bound values $\mathbf{F}^0$ and $\mathbf{M}^0$. 
    
    \item Use a fast heuristic to find a feasible solution of the OTS and obtain an upper bound $\overline{C}$ on the optimal generation cost.
    
    \item For each line $l \in \mathcal{L}$, solve the bounding problems~\eqref{eq:bounding_new} using the partially relaxed feasible region  $\mathcal{X}^k_l$.
    
    \item Solve the resulting OTS problem~\eqref{eq:ots_mip} with tighter bounds.
\end{enumerate}

\State \textbf{Output:} Optimal topology and dispatch solution.
\end{algorithmic}
\end{small}
\end{algorithm}

\section{Comparison} \label{sec:comparisson}

To evaluate the performance of the proposed methodology, we compare it against several alternative approaches. In addition to the optimization-based bound tightening procedure introduced in~\cite{pineda2024tight} (denoted as \texttt{TBT-0}), we consider three benchmark methods that solve the OTS problem without applying any bound tightening:

\begin{itemize}
    \item \texttt{MIP}: The OTS problem~\eqref{eq:ots_mip} is solved directly as a mixed-integer linear program using the original, unrefined big-M bounds proposed in \cite{moulin2010transmission}.
    
\item \texttt{BLP}: The OTS problem~\eqref{eq:OTS_NP} is solved directly as a mixed-integer bilinear program. 
All constraints (1a)–(1h) are enforced explicitly, and in particular the bilinear relation 
$f_l = x_l \tilde f_l$ in~(1b) is handled natively by the solver.

\item \texttt{IND}: The bilinear constraint~(1b) is replaced by indicator constraints linking the same variables 
$f_l$, $x_l$, and $\tilde f_l$. Specifically, if $x_l = 0$ then $f_l = 0$, and if $x_l = 1$ then 
$f_l = \tilde f_l$. This formulation avoids big-$M$ constants at the expense of increased solver effort 
\cite{klotz2013practical}.

\end{itemize}

Apart from the three benchmark approaches that do not involve bound tightening, we also compare the performance of our topological bound tightening method with an alternative strategy that leverages the full capabilities of mixed-integer linear solvers to compute tighter bounds. This alternative approach solves the bounding problems~\eqref{eq:bounding_new} by replacing the partially relaxed feasible set $\mathcal{X}^k_l$ with the fully integer-constrained set $\mathcal{X}^B$, i.e., all switching decisions are treated as binary variables. Since these bounding problems are as computationally challenging as the original OTS problem, we impose a time limit and take the best bound obtained by the solver within that time. We refer to this approach as \texttt{SBT-$t$}, where $t$ denotes the maximum solver time in milliseconds allocated to each bounding problem, and \texttt{SBT} stands for \emph{Solver Bound Tightening}.

We evaluate all discussed approaches along two main dimensions: the improvement in variable bounds and the quality and efficiency of the obtained solutions. To assess the effectiveness of the bound tightening procedures, we evaluate the relative improvement in the bounds of both the power flow variables and the dummy flow variables. Let the original bounds be denoted by $\underline{f}^0_l$, $\overline{f}^0_l$, $\underline{M}^0_l$, and $\overline{M}^0_l$. The relative reductions in bound widths are computed as:
\begin{subequations}
\begin{align}
    & \Delta F = \frac{100}{|\mathcal{L}|} \sum_{l \in \mathcal{L}} \left( 1 - \frac{\overline{f}_l - \underline{f}_l}{\overline{f}^0_l - \underline{f}^0_l} \right) \\
    & \Delta M = \frac{100}{|\mathcal{L}|} \sum_{l \in \mathcal{L}} \left( 1 - \frac{\overline{M}_l - \underline{M}_l}{\overline{M}^0_l - \underline{M}^0_l} \right)
\end{align}
\end{subequations}
These metrics capture the average percentage reduction in the bounds across all lines, providing a quantitative measure of the tightening effect.

We assess solution quality using three complementary performance metrics, all reported as relative percentages. First, the optimality gap (\texttt{gap}) measures the difference between the best known feasible solution (incumbent) and the best bound obtained by the solver. This metric reflects the solver’s ability to certify optimality. Second, the suboptimality gap (\texttt{sub}) quantifies how close the solution is to the best one found across all approaches. Specifically, it is the relative difference between the cost of the solution returned by the method and the lowest cost achieved by any approach for the same instance. Finally, the numerical discrepancy metric (\texttt{dif}) evaluates the numerical consistency of the solver-reported objective. It is computed as the relative difference between the original objective value of the reported incumbent and the cost obtained by re-solving a linear DC-OPF with the same line status variables fixed. Values of \texttt{dif} close to zero indicate high numerical reliability \cite{neumaier2004safe}.

To quantify computational burden, we report three time-related indicators: (i) the time spent solving the bound tightening problems, denoted as $T^B$; (ii) the time spent solving the resulting OTS problem, denoted as $T^O$; and (iii) the total time $T^T := T^B + T^O$. Additionally, we report the number of problem instances for which the time limit was reached, denoted as \#TL.

\section{Computational Results} \label{sec:results}

\subsection{Experimental setup}

This section presents the computational results obtained by applying the methods described in Section~\ref{sec:boundtightening} to a realistic power system. Specifically, we use the IEEE 118-bus network, which contains 186 transmission lines, as documented in~\cite{blumsack2006network}. All lines in this network are treated as switchable, as previously discussed. The 118-bus system is widely recognized as a standard benchmark in the literature on optimal transmission switching (OTS)~\cite{fisher2008optimal, kocuk2016cycle, fattahi2019bound, fuller2012fast, crozier2022feasible, hinneck2022optimal, johnson2020knearest, yang2019line, dey2022node}, as it presents a realistic and sufficiently large setting to challenge state-of-the-art optimization techniques while remaining tractable in practice.

To evaluate the performance of each method across a diverse set of scenarios, we consider 300 different OTS instances generated by varying the nodal demands. For each instance, the demand at each bus is independently sampled from a uniform distribution in the interval $[0.9\widehat{d}_n, 1.1\widehat{d}_n]$, where $\widehat{d}_n$ is the baseline demand. This sampling strategy ensures that the results reflect a broad spectrum of system operating conditions and problem complexities \cite{aguilar2025graph}. All data is available at the repository \cite{Oasys2025repository}.

All optimization problems are solved using Gurobi 10.0.3 via the \texttt{gurobipy} package version 12.0.1 in Python 3.11.4. Simulations are run on a Linux-based server equipped with an AMD EPYC processor clocking at 2.25~GHz, using a single thread and 8~GB of RAM. In all cases, the optimality gap is set to $0.01\%$ for any mixed-integer problem solved, and the time limit is set to one hour.

The bound tightening approach proposed in Section~\ref{sec:boundtightening} relies on the availability of an upper bound on the total generation cost, which must be derived from a feasible solution to the OTS problem. Although heuristic methods from the literature (e.g.,~\cite{crozier2022feasible}) could be employed for this purpose, in our study we leverage the built-in general-purpose heuristics provided by Gurobi. Specifically, we solve the OTS model with the \texttt{Heuristics} parameter set to 100\% and impose a time limit of 10 seconds. The best incumbent solution found within this time window is then used to compute the upper bound cost, denoted as~$\overline{C}$. If no feasible solution is found, the upper bound is conservatively estimated as the total system demand multiplied by the cost of the most expensive generator.

We consider five benchmark methods to evaluate the impact of bound tightening. \texttt{IND} and \texttt{BLP}, do not use any big-M values as bounds for the dummy variables, while \texttt{MIP} does require such big-M bounds to linearize the problem. Finally, we include two additional variants, \texttt{IND+} and \texttt{BLP+}, which are identical to \texttt{IND} and \texttt{BLP} but include big-M bounds on the dummy variables. In addition to these benchmarks, we evaluate three methods that solve the bounding problems including all binary variables, with a limited time budget of 25, 50, and 75 milliseconds per bounding problem. These are denoted as \texttt{SBT-25}, \texttt{SBT-50}, and \texttt{SBT-75}, and are configured to ensure a total bounding time comparable to that of the proposed method. Regarding the proposed topology-based bound tightening approach, we consider five variants with different values of the closeness parameter~$k$, resulting in methods \texttt{TBT-1}, \texttt{TBT-2}, \texttt{TBT-3}, \texttt{TBT-4}, and \texttt{TBT-5}. In these cases, the maximum time allocated to each bounding problem is limited to 5 seconds, although this limit is rarely reached in practice, as the problems typically involve only a small subset of binary variables. The computational results for all thirteen approaches are presented in Table~\ref{tab:complete_results}, where we use~$\widehat{(\cdot)}$ to denote average values over the 300 problem instances and~$\overline{(\cdot)}$ to denote maximum values.

\subsection{Results and discussion}

Focusing on the first five rows of the Table~\ref{tab:complete_results}, we compare the benchmark approaches \texttt{MIP}, \texttt{IND}, \texttt{BLP}, and their enhanced versions, \texttt{IND+} and \texttt{BLP+}. Among the original formulations, \texttt{MIP} performs noticeably better than \texttt{IND} and \texttt{BLP}, as these two approaches suffer from poorer solution quality and limited solver robustness, likely due to numerical issues arising from the absence of bounds on key variables within the indicator or nonlinear constraints. However, when such bounds are included (as in \texttt{IND+} and \texttt{BLP+}) both solution accuracy and computational efficiency improve significantly. In particular, \texttt{IND+} emerges as the best-performing benchmark overall.

\begin{table*}[]
\renewcommand{\arraystretch}{1.2}
\centering
\begin{tabular}{lcccccccccccc}
\toprule
Approach & $\widehat{\Delta F}$ & $\widehat{\Delta M}$ & $\widehat{\texttt{gap}}$ & $\overline{\texttt{gap}}$ & $\widehat{\texttt{dif}}$ & $\overline{\texttt{dif}}$ & $\widehat{\texttt{sub}}$ & $\overline{\texttt{sub}}$ & $\widehat{T^B}$(s) & $\widehat{T^O}$(s) & $\widehat{T^T}$(s) & $\#$TL \\
\midrule
\texttt{MIP} & 0.00 & 0.00 & 0.05 & 9.86 & 0.00 & 0.04 & 0.01 & 0.16 & 0.00 & 327.22 & 327.22 & 20 \\
\texttt{IND} & 0.00 & 0.00 & 1.70 & 17.78 & 1.42 & 17.29 & 0.89 & 12.88 & 0.00 & 1026.21 & 1026.21 & 60 \\
\texttt{IND+} & 0.00 & 0.00 & 0.05 & 8.54 & 0.00 & 0.15 & 0.01 & 0.23 & 0.00 & 251.98 & 251.98 & 17 \\
\texttt{BLP} & 0.00 & 0.00 & 2.20 & 14.20 & 1.35 & 9.44 & 1.18 & 11.57 & 0.00 & 2233.53 & 2233.53 & 160 \\
\texttt{BLP+} & 0.00 & 0.00 & 0.02 & 0.88 & 0.01 & 0.87 & 0.01 & 0.88 & 0.00 & 366.39 & 366.39 & 22 \\
\texttt{TBT-0} & 10.11 & 0.00 & 0.03 & 2.94 & 0.00 & 0.05 & 0.01 & 0.31 & 8.78 & 298.89 & 307.67 & 19 \\
\texttt{TBT-1} & 12.57 & 1.50 & 0.01 & 0.83 & 0.00 & 0.06 & 0.01 & 0.38 & 12.77 & 257.94 & 270.71 & 16 \\
\texttt{TBT-2} & 14.23 & 3.44 & 0.01 & 0.64 & 0.00 & 0.16 & 0.01 & 0.26 & 20.80 & 159.88 & 180.68 & 9 \\
\texttt{TBT-3} & 17.91 & 4.82 & 0.01 & 0.33 & 0.00 & 0.33 & 0.01 & 0.29 & 42.18 & 237.16 & 279.34 & 13 \\
\texttt{TBT-4} & 20.50 & 5.89 & 0.01 & 0.16 & 0.00 & 0.05 & 0.01 & 0.10 & 145.33 & 133.67 & 279.00 & 8 \\
\texttt{TBT-5} & 21.08 & 10.00 & 0.01 & 1.63 & 0.00 & 0.03 & 0.00 & 0.08 & 384.45 & 136.56 & 521.01 & 6 \\
\texttt{SBT-25} & 11.81 & 2.09 & 0.01 & 0.22 & 0.00 & 0.20 & 0.00 & 0.04 & 18.29 & 216.83 & 235.12 & 11 \\
\texttt{SBT-50} & 12.06 & 2.61 & 0.02 & 0.46 & 0.00 & 0.14 & 0.01 & 0.13 & 36.56 & 295.47 & 332.03 & 19 \\
\texttt{SBT-75} & 12.49 & 2.89 & 0.01 & 0.52 & 0.00 & 0.05 & 0.01 & 0.23 & 54.54 & 231.85 & 286.40 & 14 \\
\bottomrule
\end{tabular}
\caption{Computational results for all instances.}
\label{tab:complete_results}
\end{table*}

Next, we analyze the results of the proposed topology-based bound tightening (TBT) methodology for different values of the closeness parameter~$k$. As expected, increasing~$k$ leads to tighter average bound reductions, since the corresponding bounding problems become less relaxed ($\widehat{\Delta F}, \widehat{\Delta M}$). However, this comes at the expense of a higher average computational time required to solve these bounding problems ($\widehat{T^B}$). The growth of the closeness set with $k$ helps to understand this trade-off: the average number of lines in the set $|\mathcal{X}^k|$ is 6.38, 20.01, 39.82, 62.42, and 86.34 for $k=1$ to $5$, respectively. It is also noteworthy that the approach proposed in~\cite{pineda2024tight}, denoted as \texttt{TBT-0}, is unable to tighten the big-M bounds when all lines are switchable ($\Delta M = 0.00$).

Second, although increasing~$k$ consistently produces tighter bounds, this does not always translate into a shorter solution time for the resulting OTS instances ($\widehat{T^O}$). For instance, the average time to solve the OTS is higher for $k=3$ than for $k=2$. While this may seem counterintuitive, it can be explained by the inherent performance variability in solving mixed-integer programs, where even minor changes (such as permuting two rows or columns of the constraint matrix) can significantly affect the solver's behavior \cite{lodi2013performance}. In our case, although tighter bounds are generally beneficial, they can influence branching decisions, cut generation, or other internal heuristics of the solver, potentially triggering snowballing effects that result in longer solution times \cite{le2015important, wojtaszek2010faster}.

Despite these uncontrollable factors, we observe that all proposed \texttt{TBT} variants yield better overall performance than the baseline \texttt{MIP} approach. Specifically, they achieve lower maximum and average optimality gaps, shorter total computational times, and fewer unsolved instances within the one-hour time limit. Focusing on the average total time ($\widehat{T^T}$), \texttt{TBT-0} reduces it by 6\% and solves one additional instance compared to \texttt{MIP}, while the proposed method \texttt{TBT-2} achieves a 45\% reduction and cuts the number of unsolved problems by half. Although \texttt{TBT-5} solves the highest number of instances, its average time increases due to the additional effort required for solving more complex bounding problems. When comparing with the enhanced formulation \texttt{IND+}, we observe that all \texttt{TBT} methods further reduce the number of unsolved instances and the maximum gap. However, the average total time is only improved by \texttt{TBT-2}, which reduces it by 28\% relative to \texttt{IND+}. Finally, the average metrics for solution suboptimality (\texttt{sub}) and numerical consistency (\texttt{dif}) of the proposed methods remain comparable to those of both \texttt{MIP} and \texttt{IND+}, confirming that the gains in performance are not achieved at the expense of solution quality.

Next, we analyze the results for the solver bound tightening approaches (SBT), where the bounding problems include all binary variables but are subject to time limits of 25, 50, or 75 milliseconds. First, we observe that the performance of \texttt{SBT-50} closely resembles that of the baseline \texttt{MIP} method in terms of average total solution time and the number of unsolved instances. Furthermore, although \texttt{SBT-75} dedicates more time to solving the bounding problems than \texttt{TBT-3}, the resulting bound reductions are notably smaller. This suggests that selectively relaxing a subset of switching variables based on topological proximity is a more effective strategy than attempting to solve the full binary problem under stringent time constraints. 

\begin{figure}
    \centering
\begin{tikzpicture}
\begin{axis}[
    width=8cm, height=5.5cm,
    xlabel={Total time (s)},
    ylabel style={yshift=-1.5em},
    ylabel={Number of instances solved},
    xmin=0, xmax=3671.93,
    ymin=0, ymax=101.98,
    legend style={at={(0.99,0.01)},anchor=south east, font=\footnotesize},
    grid=major,
]
\addplot+[no markers, thick, color=blue, const plot] coordinates {%
(3.46,1) (5.15,2) (8.13,3) (8.51,4) (8.61,5) (9.20,6) (9.33,7) (10.29,8) (10.47,9) (10.48,10) (10.70,11) (11.11,12) (12.57,13) (12.97,14) (13.42,15) (14.68,16) (16.77,17) (17.32,18) (17.68,19) (17.97,20) (18.54,21) (18.87,22) (19.33,23) (19.36,24) (20.95,25) (21.48,26) (23.46,27) (23.63,28) (23.74,29) (24.73,30) (24.81,31) (25.26,32) (25.87,33) (26.14,34) (26.47,35) (26.48,36) (28.10,37) (28.28,38) (28.78,39) (29.04,40) (29.05,41) (33.95,42) (34.46,43) (34.59,44) (35.43,45) (36.28,46) (39.11,47) (39.22,48) (39.75,49) (41.14,50) (42.29,51) (47.05,52) (49.54,53) (56.13,54) (58.53,55) (59.84,56) (63.44,57) (68.71,58) (71.95,59) (72.64,60) (77.48,61) (83.29,62) (92.15,63) (102.28,64) (102.40,65) (106.33,66) (106.62,67) (109.41,68) (111.92,69) (120.99,70) (121.30,71) (148.49,72) (185.85,73) (288.42,74) (330.56,75) (381.68,76) (484.35,77) (760.98,78) (989.40,79) (1134.86,80) (1227.62,81) (1833.23,82) (2308.61,83) (2509.39,84) (3600.00,85) (3600.00,86) (3600.00,87) (3600.00,88) (3600.00,89) (3600.00,90) (3600.00,91) (3600.00,92) (3600.00,93) (3600.00,94) (3600.00,95) (3600.00,96) (3600.00,97) (3600.00,98) (3600.00,99) (3600.00,100) };
\addlegendentry{\texttt{MIP}}
\addplot+[no markers, thick, color=red, const plot] coordinates {%
(5.15,1) (6.20,2) (6.91,3) (8.05,4) (8.50,5) (8.81,6) (9.77,7) (10.94,8) (11.36,9) (12.24,10) (12.29,11) (12.31,12) (12.36,13) (13.48,14) (14.21,15) (14.40,16) (14.44,17) (14.71,18) (15.94,19) (15.99,20) (16.04,21) (16.46,22) (16.72,23) (17.32,24) (17.51,25) (17.68,26) (17.89,27) (18.44,28) (19.21,29) (20.58,30) (21.01,31) (22.48,32) (22.63,33) (22.92,34) (23.38,35) (24.13,36) (24.63,37) (24.64,38) (24.64,39) (24.99,40) (25.37,41) (26.72,42) (27.15,43) (27.26,44) (29.10,45) (29.57,46) (29.96,47) (31.24,48) (32.02,49) (33.13,50) (33.43,51) (33.45,52) (34.06,53) (36.01,54) (36.74,55) (38.09,56) (39.03,57) (40.72,58) (40.81,59) (43.61,60) (43.87,61) (46.65,62) (47.34,63) (47.51,64) (47.59,65) (49.45,66) (51.10,67) (51.79,68) (53.92,69) (56.05,70) (62.04,71) (84.45,72) (91.56,73) (106.95,74) (119.67,75) (168.55,76) (233.41,77) (238.33,78) (361.42,79) (415.91,80) (473.58,81) (1016.28,82) (1248.17,83) (3600.00,84) (3600.00,85) (3600.00,86) (3600.00,87) (3600.00,88) (3600.00,89) (3600.00,90) (3600.00,91) (3600.00,92) (3600.00,93) (3600.00,94) (3600.00,95) (3600.00,96) (3600.00,97) (3600.00,98) (3600.00,99) (3600.00,100) };
\addlegendentry{\texttt{IND+}}
\addplot+[no markers, thick, color=green!60!black, const plot] coordinates {%
(15.16,1) (15.55,2) (15.90,3) (16.56,4) (17.50,5) (17.85,6) (18.64,7) (18.68,8) (19.15,9) (19.45,10) (20.42,11) (21.21,12) (21.22,13) (23.20,14) (23.87,15) (24.86,16) (24.95,17) (25.15,18) (25.41,19) (25.81,20) (26.19,21) (26.81,22) (26.94,23) (27.27,24) (27.53,25) (27.57,26) (27.75,27) (28.96,28) (29.04,29) (29.06,30) (29.28,31) (29.34,32) (33.76,33) (33.85,34) (34.18,35) (35.55,36) (36.63,37) (36.99,38) (37.37,39) (39.02,40) (39.17,41) (40.10,42) (40.62,43) (40.85,44) (40.94,45) (43.56,46) (44.83,47) (46.62,48) (47.81,49) (48.61,50) (49.60,51) (49.94,52) (51.17,53) (53.49,54) (53.86,55) (56.31,56) (56.72,57) (56.87,58) (57.40,59) (60.33,60) (60.89,61) (61.12,62) (61.35,63) (61.84,64) (65.55,65) (67.09,66) (67.63,67) (81.15,68) (87.57,69) (90.37,70) (95.78,71) (105.49,72) (113.61,73) (114.96,74) (116.23,75) (129.61,76) (187.98,77) (199.47,78) (244.63,79) (245.30,80) (627.97,81) (2481.36,82) (2684.54,83) (3459.27,84) (3600.00,85) (3600.00,86) (3600.00,87) (3600.00,88) (3600.00,89) (3600.00,90) (3600.00,91) (3600.00,92) (3600.00,93) (3600.00,94) (3600.00,95) (3600.00,96) (3600.00,97) (3600.00,98) (3600.00,99) (3600.00,100) };
\addlegendentry{\texttt{TBT-0}}
\addplot+[no markers, thick, color=orange, const plot] coordinates {%
(24.89,1) (25.22,2) (26.37,3) (26.98,4) (28.10,5) (28.80,6) (29.41,7) (29.48,8) (29.52,9) (29.52,10) (29.65,11) (29.85,12) (31.42,13) (31.66,14) (31.76,15) (31.76,16) (31.82,17) (32.23,18) (32.47,19) (32.91,20) (33.13,21) (33.72,22) (34.01,23) (34.34,24) (34.89,25) (35.20,26) (35.44,27) (35.52,28) (35.58,29) (35.69,30) (36.19,31) (38.31,32) (38.95,33) (39.65,34) (40.31,35) (40.45,36) (40.47,37) (40.74,38) (42.99,39) (44.52,40) (45.93,41) (46.00,42) (46.08,43) (46.15,44) (46.31,45) (46.41,46) (46.49,47) (47.47,48) (48.23,49) (48.29,50) (48.50,51) (48.61,52) (48.96,53) (50.49,54) (52.67,55) (53.20,56) (53.42,57) (55.14,58) (55.39,59) (55.55,60) (57.50,61) (57.79,62) (57.88,63) (58.09,64) (60.25,65) (61.24,66) (62.82,67) (63.35,68) (65.53,69) (67.69,70) (69.10,71) (69.29,72) (69.49,73) (70.61,74) (71.23,75) (71.55,76) (75.36,77) (76.87,78) (79.77,79) (90.84,80) (90.95,81) (92.00,82) (100.65,83) (109.07,84) (114.26,85) (132.51,86) (221.72,87) (223.94,88) (297.45,89) (326.07,90) (493.09,91) (512.04,92) (513.86,93) (2632.40,94) (3600.00,95) (3600.00,96) (3600.00,97) (3600.00,98) (3600.00,99) (3600.00,100) };
\addlegendentry{\texttt{TBT-2}}
\addplot+[no markers, thick, color=purple, const plot] coordinates {%
(135.36,1) (138.68,2) (146.32,3) (148.00,4) (150.34,5) (150.74,6) (151.27,7) (152.52,8) (152.86,9) (152.91,10) (153.46,11) (153.76,12) (153.78,13) (154.22,14) (154.23,15) (154.62,16) (155.14,17) (155.69,18) (155.72,19) (156.02,20) (156.08,21) (157.33,22) (157.49,23) (157.61,24) (157.66,25) (157.70,26) (157.76,27) (157.79,28) (157.80,29) (158.77,30) (159.20,31) (159.38,32) (159.69,33) (159.77,34) (159.87,35) (159.90,36) (160.71,37) (160.96,38) (163.04,39) (164.15,40) (164.46,41) (165.09,42) (166.64,43) (166.84,44) (166.96,45) (167.46,46) (168.31,47) (168.42,48) (168.80,49) (169.13,50) (169.60,51) (170.31,52) (170.55,53) (172.37,54) (172.54,55) (173.08,56) (173.28,57) (174.03,58) (174.97,59) (175.04,60) (175.74,61) (175.89,62) (176.77,63) (177.13,64) (177.47,65) (178.13,66) (182.87,67) (183.70,68) (184.32,69) (184.63,70) (184.74,71) (186.47,72) (187.69,73) (192.70,74) (195.49,75) (196.11,76) (196.30,77) (200.74,78) (211.12,79) (213.03,80) (220.41,81) (221.93,82) (223.29,83) (241.63,84) (281.30,85) (366.64,86) (379.33,87) (393.20,88) (475.49,89) (597.27,90) (606.68,91) (652.13,92) (3600.00,93) (3600.00,94) (3600.00,95) (3600.00,96) (3600.00,97) (3600.00,98) (3600.00,99) (3600.00,100) };
\addlegendentry{\texttt{TBT-4}}
\addplot+[no markers, thick, color=cyan, const plot, solid] coordinates {%
(21.20,1) (23.28,2) (23.38,3) (24.52,4) (25.61,5) (26.48,6) (26.63,7) (26.88,8) (26.94,9) (27.12,10) (27.30,11) (27.43,12) (27.61,13) (27.68,14) (28.33,15) (28.87,16) (29.53,17) (29.56,18) (29.83,19) (30.30,20) (30.58,21) (30.73,22) (32.51,23) (32.59,24) (33.11,25) (33.16,26) (33.39,27) (33.49,28) (33.95,29) (34.14,30) (34.59,31) (35.05,32) (35.11,33) (35.20,34) (35.29,35) (35.44,36) (35.81,37) (37.17,38) (39.69,39) (39.83,40) (40.50,41) (41.35,42) (42.11,43) (42.64,44) (44.26,45) (44.36,46) (44.77,47) (46.73,48) (47.21,49) (47.56,50) (50.96,51) (53.58,52) (54.26,53) (54.28,54) (54.65,55) (56.15,56) (58.15,57) (58.41,58) (58.59,59) (60.35,60) (62.95,61) (70.38,62) (72.11,63) (73.54,64) (87.93,65) (91.45,66) (91.66,67) (102.28,68) (103.26,69) (130.56,70) (132.86,71) (135.14,72) (167.76,73) (172.16,74) (193.55,75) (207.74,76) (211.90,77) (222.04,78) (227.00,79) (232.84,80) (247.51,81) (283.52,82) (289.47,83) (294.52,84) (322.57,85) (443.00,86) (453.73,87) (1243.64,88) (1594.05,89) (2877.85,90) (3600.00,91) (3600.00,92) (3600.00,93) (3600.00,94) (3600.00,95) (3600.00,96) (3600.00,97) (3600.00,98) (3600.00,99) (3600.00,100) };
\addlegendentry{\texttt{SBT-25}}
\end{axis}
\end{tikzpicture}
    \caption{Performance profiles for hard instances.}
    \label{fig:hard_instances}
\end{figure}
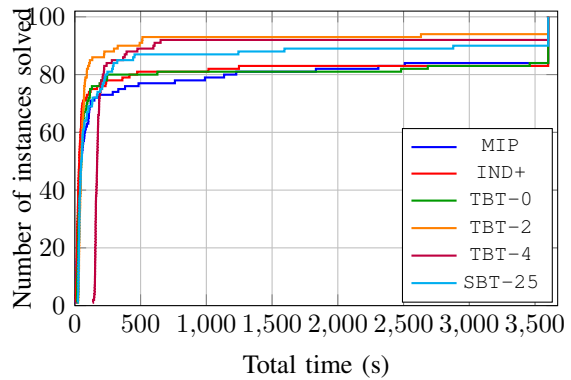

We continue the analysis of the computational results by selecting a subset of 100 \emph{hard} instances and 100 \emph{easy} instances from the original set of 300. To this end, we compute for each instance the average solving time across all approaches considered. Instances with the highest average times are classified as hard, while those with the lowest are classified as easy. In this part of the analysis, we reduce the number of approaches to a subset including the benchmark methods (\texttt{MIP}, \texttt{IND+} and \texttt{TBT-0}), the two variants of the proposed approach that offer the best trade-off between average total time and number of unsolved instances (\texttt{TBT-2} and \texttt{TBT-4}), and the variant of the solver-based bound tightening method with the lowest average total time and fewest unsolved problems (\texttt{SBT-25}). Fig.~\ref{fig:easy_instances} displays, for these selected approaches, the number of easy instances solved within a given time limit, while Fig.~\ref{fig:hard_instances} presents the corresponding results for the hard instances.

\begin{figure}
    \centering
\begin{tikzpicture}
\begin{axis}[
    width=8cm, height=5.5cm,
    xlabel={Total time (s)},
    ylabel style={yshift=-1.5em},
    ylabel={Number of instances solved},
    xmin=0, xmax=1200,
    ymin=0, ymax=101.98,
    legend style={at={(0.99,0.01)},anchor=south east, font=\footnotesize},
    grid=major,
]
\addplot+[no markers, thick, color=blue, const plot] coordinates {%
(3.13,1) (4.11,2) (4.31,3) (4.89,4) (5.05,5) (6.07,6) (6.78,7) (6.94,8) (7.09,9) (7.27,10) (7.41,11) (7.91,12) (8.00,13) (8.43,14) (8.75,15) (8.75,16) (8.80,17) (8.94,18) (9.08,19) (9.41,20) (10.25,21) (10.57,22) (10.65,23) (10.76,24) (11.11,25) (11.56,26) (11.60,27) (11.82,28) (12.19,29) (12.22,30) (12.66,31) (12.74,32) (12.99,33) (13.85,34) (13.93,35) (14.20,36) (14.24,37) (14.52,38) (15.61,39) (15.96,40) (16.10,41) (16.23,42) (16.61,43) (16.82,44) (17.86,45) (18.19,46) (18.52,47) (18.61,48) (18.77,49) (18.82,50) (18.98,51) (19.09,52) (19.56,53) (19.65,54) (21.28,55) (21.34,56) (21.35,57) (22.16,58) (23.62,59) (23.78,60) (23.80,61) (23.97,62) (24.01,63) (25.06,64) (25.08,65) (25.14,66) (25.31,67) (25.75,68) (29.08,69) (29.41,70) (29.96,71) (30.23,72) (31.24,73) (32.19,74) (33.12,75) (33.41,76) (35.71,77) (36.28,78) (37.31,79) (38.10,80) (38.16,81) (38.81,82) (38.84,83) (39.16,84) (39.23,85) (39.30,86) (39.46,87) (42.11,88) (47.94,89) (59.47,90) (64.21,91) (64.66,92) (72.65,93) (88.55,94) (126.03,95) (145.46,96) (168.34,97) (589.22,98) (1010.14,99) (3600.00,100) };
\addlegendentry{\texttt{MIP}}
\addplot+[no markers, thick, color=red, const plot] coordinates {%
(2.89,1) (3.00,2) (3.47,3) (3.69,4) (3.87,5) (5.06,6) (5.84,7) (5.91,8) (6.56,9) (6.59,10) (6.64,11) (6.74,12) (6.97,13) (7.10,14) (7.48,15) (7.56,16) (7.62,17) (7.69,18) (8.43,19) (8.44,20) (8.74,21) (8.92,22) (9.64,23) (9.85,24) (9.91,25) (10.08,26) (10.32,27) (10.38,28) (10.51,29) (10.56,30) (10.71,31) (11.56,32) (11.80,33) (12.19,34) (12.20,35) (12.38,36) (12.92,37) (13.37,38) (13.46,39) (13.50,40) (13.62,41) (14.02,42) (14.05,43) (14.38,44) (14.69,45) (14.74,46) (16.24,47) (16.29,48) (16.66,49) (17.46,50) (18.09,51) (18.79,52) (19.08,53) (19.55,54) (19.72,55) (20.36,56) (20.46,57) (20.71,58) (20.78,59) (21.15,60) (21.49,61) (21.85,62) (21.97,63) (22.51,64) (22.70,65) (23.73,66) (24.76,67) (25.40,68) (25.48,69) (25.75,70) (26.65,71) (28.42,72) (28.73,73) (29.18,74) (29.48,75) (29.78,76) (29.88,77) (30.70,78) (30.76,79) (31.37,80) (31.53,81) (31.95,82) (33.23,83) (34.10,84) (34.14,85) (35.10,86) (35.95,87) (36.47,88) (41.19,89) (42.05,90) (43.08,91) (45.32,92) (50.74,93) (54.20,94) (55.19,95) (70.04,96) (73.13,97) (132.29,98) (149.03,99) (451.32,100) };
\addlegendentry{\texttt{IND+}}
\addplot+[no markers, thick, color=green!60!black, const plot] coordinates {%
(11.40,1) (11.80,2) (13.66,3) (13.93,4) (14.05,5) (14.15,6) (14.25,7) (14.82,8) (15.24,9) (15.25,10) (15.49,11) (15.72,12) (15.74,13) (15.79,14) (15.80,15) (15.83,16) (15.97,17) (16.13,18) (16.25,19) (16.33,20) (17.64,21) (17.67,22) (17.68,23) (18.20,24) (18.47,25) (18.81,26) (19.14,27) (19.43,28) (19.69,29) (19.93,30) (20.07,31) (20.27,32) (20.37,33) (20.38,34) (20.52,35) (20.69,36) (20.95,37) (21.05,38) (21.17,39) (21.27,40) (21.59,41) (21.65,42) (22.19,43) (22.28,44) (22.28,45) (22.73,46) (23.16,47) (23.74,48) (24.16,49) (24.76,50) (24.86,51) (24.89,52) (25.39,53) (25.40,54) (25.99,55) (26.18,56) (26.19,57) (26.42,58) (26.55,59) (26.64,60) (26.96,61) (27.04,62) (27.18,63) (27.19,64) (28.56,65) (28.95,66) (29.13,67) (30.71,68) (31.45,69) (32.01,70) (32.18,71) (33.28,72) (33.48,73) (33.84,74) (34.33,75) (34.47,76) (35.26,77) (36.96,78) (37.25,79) (37.91,80) (38.74,81) (39.48,82) (40.58,83) (41.82,84) (45.90,85) (49.46,86) (50.62,87) (52.15,88) (52.36,89) (52.59,90) (55.96,91) (57.23,92) (73.36,93) (79.06,94) (127.56,95) (190.12,96) (310.88,97) (490.46,98) (531.56,99) (551.23,100) };
\addlegendentry{\texttt{TBT-0}}
\addplot+[no markers, thick, color=orange, const plot] coordinates {%
(24.05,1) (24.93,2) (25.38,3) (26.26,4) (26.44,5) (26.58,6) (26.66,7) (26.78,8) (26.85,9) (27.19,10) (27.99,11) (28.15,12) (28.17,13) (28.20,14) (28.27,15) (28.80,16) (30.02,17) (30.21,18) (30.68,19) (30.73,20) (30.84,21) (31.89,22) (32.03,23) (32.15,24) (32.26,25) (32.50,26) (32.50,27) (32.96,28) (33.00,29) (33.18,30) (33.23,31) (33.25,32) (33.65,33) (33.92,34) (33.95,35) (34.05,36) (34.45,37) (34.88,38) (34.98,39) (35.58,40) (35.74,41) (35.74,42) (35.96,43) (36.10,44) (36.19,45) (36.52,46) (36.77,47) (36.77,48) (37.04,49) (37.18,50) (37.32,51) (37.46,52) (37.58,53) (37.61,54) (38.64,55) (38.67,56) (39.08,57) (39.40,58) (40.02,59) (41.36,60) (42.08,61) (42.19,62) (42.66,63) (43.19,64) (43.24,65) (43.36,66) (44.93,67) (45.12,68) (45.68,69) (46.48,70) (47.59,71) (48.21,72) (50.09,73) (51.71,74) (51.76,75) (52.31,76) (53.07,77) (53.21,78) (53.86,79) (56.42,80) (59.38,81) (59.55,82) (59.86,83) (60.34,84) (62.51,85) (66.03,86) (66.18,87) (69.40,88) (71.01,89) (78.40,90) (101.76,91) (102.91,92) (122.78,93) (147.68,94) (160.32,95) (178.38,96) (214.12,97) (270.89,98) (686.87,99) (1070.60,100) };
\addlegendentry{\texttt{TBT-2}}
\addplot+[no markers, thick, color=purple, const plot] coordinates {%
(125.15,1) (126.70,2) (128.79,3) (133.59,4) (139.96,5) (141.52,6) (141.86,7) (142.94,8) (143.81,9) (144.10,10) (145.95,11) (148.15,12) (148.75,13) (149.85,14) (151.09,15) (151.48,16) (151.70,17) (151.97,18) (152.43,19) (152.67,20) (155.36,21) (155.42,22) (155.50,23) (155.53,24) (155.82,25) (156.41,26) (156.45,27) (157.59,28) (157.99,29) (158.69,30) (159.01,31) (159.88,32) (160.41,33) (160.60,34) (160.88,35) (161.19,36) (161.26,37) (161.43,38) (161.58,39) (161.92,40) (162.99,41) (163.58,42) (164.40,43) (164.42,44) (164.91,45) (165.41,46) (166.33,47) (166.42,48) (166.48,49) (166.82,50) (167.72,51) (167.98,52) (168.20,53) (168.35,54) (168.80,55) (168.86,56) (168.95,57) (169.61,58) (170.53,59) (170.60,60) (170.64,61) (171.06,62) (172.04,63) (172.05,64) (173.98,65) (174.20,66) (174.80,67) (177.90,68) (178.37,69) (181.11,70) (181.91,71) (182.62,72) (182.63,73) (183.67,74) (183.87,75) (186.20,76) (187.02,77) (190.45,78) (191.67,79) (192.00,80) (192.94,81) (193.96,82) (194.13,83) (195.27,84) (199.96,85) (203.60,86) (204.23,87) (204.39,88) (204.95,89) (206.32,90) (208.70,91) (209.49,92) (210.09,93) (216.33,94) (217.02,95) (229.37,96) (232.12,97) (250.35,98) (253.42,99) (314.93,100) };
\addlegendentry{\texttt{TBT-4}}
\addplot+[no markers, thick, color=cyan, const plot, solid] coordinates {%
(20.76,1) (20.85,2) (21.89,3) (22.11,4) (22.52,5) (22.67,6) (23.26,7) (23.26,8) (24.56,9) (24.66,10) (24.69,11) (25.57,12) (25.63,13) (25.66,14) (25.68,15) (25.78,16) (25.81,17) (25.86,18) (26.22,19) (26.23,20) (26.26,21) (26.97,22) (27.76,23) (27.80,24) (27.82,25) (28.08,26) (28.19,27) (28.68,28) (28.97,29) (29.06,30) (29.44,31) (29.60,32) (29.60,33) (29.73,34) (29.89,35) (30.10,36) (30.26,37) (30.30,38) (30.66,39) (30.75,40) (30.85,41) (31.22,42) (31.28,43) (31.70,44) (32.20,45) (32.64,46) (32.91,47) (33.54,48) (33.73,49) (34.01,50) (34.88,51) (34.97,52) (35.13,53) (35.74,54) (36.27,55) (36.44,56) (36.65,57) (37.93,58) (38.22,59) (38.44,60) (39.76,61) (40.00,62) (40.36,63) (40.60,64) (40.74,65) (41.27,66) (42.43,67) (42.60,68) (43.22,69) (43.66,70) (45.31,71) (46.86,72) (48.00,73) (48.51,74) (49.00,75) (49.13,76) (49.97,77) (49.98,78) (50.99,79) (51.10,80) (51.44,81) (52.71,82) (52.88,83) (53.63,84) (56.35,85) (60.10,86) (60.69,87) (62.37,88) (68.78,89) (69.77,90) (71.51,91) (73.70,92) (78.04,93) (90.74,94) (122.65,95) (155.37,96) (236.79,97) (276.54,98) (367.04,99) (1167.46,100) };
\addlegendentry{\texttt{SBT-25}}
\end{axis}
\end{tikzpicture}
    \caption{Performance profiles for easy instances.}
    \label{fig:easy_instances}
\end{figure}
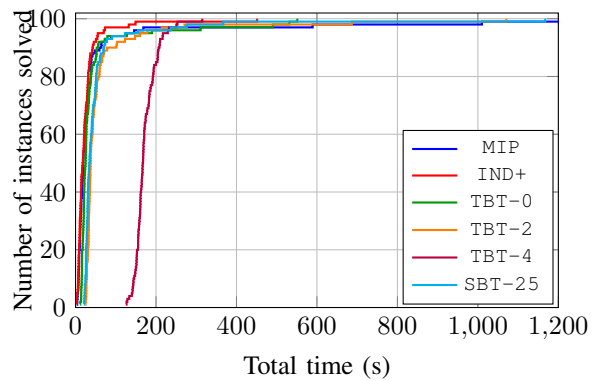

First, by comparing Fig.\ref{fig:easy_instances} and Fig.\ref{fig:hard_instances}, we observe that all methods except \texttt{MIP} are able to solve the 100 \emph{easy} instances within the 1200 seconds. In contrast, none of the methods is able to solve all 100 \emph{hard} instances within the one-hour limit, highlighting the increased complexity of this subset. Second, for the easy instances, all approaches exhibit similar performance, with only minor differences, except for \texttt{TBT-4}, which incurs a noticeable overhead due to the time spent solving the bounding problems. In contrast, for the hard instances, the benchmark approaches \texttt{MIP}, \texttt{IND+} and \texttt{TBT-0} perform significantly worse than the other methods. Their similar performance suggests that solving bounding problems with all binary variables relaxed provides little benefit when all lines are switchable. Among the remaining methods, \texttt{TBT-2} outperforms both \texttt{TBT-4} and \texttt{SBT-25}. This is also confirmed by the summarized computational results in Table~\ref{tab:hard_instances} provides further insights for the subset of hard instances, where the \texttt{TBT-2} approach achieves the lowest number of unsolved instances and reduces the total computational time by 57\% with respect to  \texttt{MIP} and by 54\% relative to \texttt{IND+}.

\begin{table}[h]
\centering
\renewcommand{\arraystretch}{1.2}
\begin{tabular}{lcccc}
\toprule
Approach & $\widehat{T^B}$ (s) & $\widehat{T^O}$ (s) & $\widehat{T^T}$ (s) & \#NL \\
\midrule
\texttt{MIP} & 0.00 & 729.83 & 729.83 & 16 \\
\texttt{IND+} & 0.00 & 676.24 & 676.24 & 17 \\
\texttt{TBT-0} & 8.78 & 702.45 & 711.23 & 16 \\
\texttt{TBT-2} & 20.83 & 290.66 & 311.49 & 6 \\
\texttt{TBT-4} & 145.49 & 322.78 & 468.27 & 8 \\
\texttt{SBT-25} & 18.29 & 472.40 & 490.69 & 10 \\
\bottomrule
\end{tabular}
\caption{Summary of computation results for hard instances.} \label{tab:hard_instances}
\end{table}

In summary, \texttt{TBT-2} performs comparably to other approaches on easy instances and significantly better on hard ones, striking an effective balance between solution quality and computational efficiency. By preserving the binary switching decisions for lines within a closeness level of 2 to the target line, this configuration maintains the enforcement of Kirchhoff’s laws in the most relevant part of the network. This selective integrality captures key electrical interactions that lead to tighter and more informative bounds, without making the bounding problems intractable. While increasing the closeness level might further improve bound strength, the added complexity outweighs the benefits. Thus \texttt{TBT-2} represents a sweet spot that retains essential network structure while ensuring computational tractability.

\section{Conclusions} \label{sec:conclusions}

Topology optimization offers valuable flexibility for power systems but is challenged by combinatorial complexity and numerical instability from large constants. Existing optimization-based bound tightening methods lose effectiveness when all lines are switchable due to overly loose relaxations. We propose a topology-aware bound tightening technique that selectively preserves binary switching decisions for lines electrically close to the one being bounded, enforcing local physical constraints and producing stronger bounds. While keeping too many binaries increases computational complexity, retaining binaries within a closeness level of two strikes a clear sweet spot, balancing bound tightness and computational effort. Computational tests on the IEEE 118-bus system show this approach halves unsolved instances within an hour and reduces average solution times by 45\%, reaching 57\% for the hardest cases, demonstrating improved scalability and efficiency for topology optimization in power systems.

\bibliographystyle{IEEEtran}
\bibliography{references}

\begin{thebibliography}{10}
\providecommand{\url}[1]{#1}
\csname url@samestyle\endcsname
\providecommand{\newblock}{\relax}
\providecommand{\bibinfo}[2]{#2}
\providecommand{\BIBentrySTDinterwordspacing}{\spaceskip=0pt\relax}
\providecommand{\BIBentryALTinterwordstretchfactor}{4}
\providecommand{\BIBentryALTinterwordspacing}{\spaceskip=\fontdimen2\font plus
\BIBentryALTinterwordstretchfactor\fontdimen3\font minus \fontdimen4\font\relax}
\providecommand{\BIBforeignlanguage}[2]{{%
\expandafter\ifx\csname l@#1\endcsname\relax
\typeout{** WARNING: IEEEtran.bst: No hyphenation pattern has been}%
\typeout{** loaded for the language `#1'. Using the pattern for}%
\typeout{** the default language instead.}%
\else
\language=\csname l@#1\endcsname
\fi
#2}}
\providecommand{\BIBdecl}{\relax}
\BIBdecl

\bibitem{o2005dispatchable}
R.~P. O'Neill, R.~Baldick, U.~Helman, M.~H. Rothkopf, and W.~Stewart, ``Dispatchable transmission in rto markets,'' \emph{IEEE Transactions on Power Systems}, vol.~20, no.~1, pp. 171--179, 2005.

\bibitem{fisher2008optimal}
E.~B. Fisher, R.~P. O'Neill, and M.~C. Ferris, ``Optimal transmission switching,'' \emph{IEEE Transactions on Power Systems}, vol.~23, no.~3, pp. 1346--1355, 2008.

\bibitem{pineda2024tight}
S.~Pineda, J.~M. Morales, {\'A}.~Porras, and C.~Dom{\'\i}nguez, ``Tight big-ms for optimal transmission switching,'' \emph{Electric Power Systems Research}, vol. 234, p. 110620, 2024.

\bibitem{numan2023role}
M.~Numan, M.~F. Abbas, M.~Yousif, S.~S. Ghoneim, A.~Mohammad, and A.~Noorwali, ``The role of optimal transmission switching in enhancing grid flexibility: A review,'' \emph{IEEE Access}, 2023.

\bibitem{fattahi2019bound}
S.~Fattahi, J.~Lavaei, and A.~Atamtürk, ``A bound strengthening method for optimal transmission switching in power systems,'' \emph{IEEE Transactions on Power Systems}, vol.~34, no.~1, pp. 280--291, 2019.

\bibitem{moulin2010transmission}
L.~S. Moulin, M.~Poss, and C.~Sagastiz{\'a}bal, ``Transmission expansion planning with re-design,'' \emph{Energy systems}, vol.~1, pp. 113--139, 2010.

\bibitem{flores2020alternative}
M.~Flores, L.~H. Macedo, and R.~Romero, ``Alternative mathematical models for the optimal transmission switching problem,'' \emph{IEEE Systems Journal}, vol.~15, no.~1, pp. 1245--1255, 2020.

\bibitem{crozier2022feasible}
C.~Crozier, K.~Baker, and B.~Toomey, ``Feasible region-based heuristics for optimal transmission switching,'' \emph{Sustainable Energy, Grids and Networks}, vol.~30, p. 100628, 2022.

\bibitem{soroush2013accuracies}
M.~Soroush and J.~D. Fuller, ``Accuracies of optimal transmission switching heuristics based on dcopf and acopf,'' \emph{IEEE Transactions on Power Systems}, vol.~29, no.~2, pp. 924--932, 2013.

\bibitem{barrows2014correcting}
C.~Barrows, S.~Blumsack, and P.~Hines, ``Correcting optimal transmission switching for ac power flows,'' in \emph{2014 47th Hawaii International Conference on System Sciences}.\hskip 1em plus 0.5em minus 0.4em\relax IEEE, 2014, pp. 2374--2379.

\bibitem{kocuk2017new}
B.~Kocuk, S.~S. Dey, and X.~A. Sun, ``New formulation and strong misocp relaxations for ac optimal transmission switching problem,'' \emph{IEEE Transactions on Power Systems}, vol.~32, no.~6, pp. 4161--4170, 2017.

\bibitem{kocuk2016cycle}
B.~Kocuk, H.~Jeon, S.~S. Dey, J.~Linderoth, J.~Luedtke, and X.~A. Sun, ``A cycle-based formulation and valid inequalities for {DC} power transmission problems with switching,'' \emph{Operations Research}, vol.~64, no.~4, pp. 922--938, 2016.

\bibitem{liu2012heuristic}
C.~Liu, J.~Wang, and J.~Ostrowski, ``Heuristic prescreening switchable branches in optimal transmission switching,'' \emph{IEEE Transactions on Power Systems}, vol.~27, no.~4, pp. 2289--2290, 2012.

\bibitem{barrows2012computationally}
C.~Barrows, S.~Blumsack, and R.~Bent, ``Computationally efficient optimal transmission switching: Solution space reduction,'' in \emph{2012 IEEE Power and Energy Society General Meeting}.\hskip 1em plus 0.5em minus 0.4em\relax IEEE, 2012, pp. 1--8.

\bibitem{fuller2012fast}
J.~D. Fuller, R.~Ramasra, and A.~Cha, ``Fast heuristics for transmission-line switching,'' \emph{IEEE Transactions on Power Systems}, vol.~27, no.~3, pp. 1377--1386, 2012.

\bibitem{hinneck2022optimal}
A.~Hinneck and D.~Pozo, ``Optimal transmission switching: improving exact algorithms by parallel incumbent solution generation,'' \emph{IEEE Transactions on Power Systems}, 2022.

\bibitem{johnson2020knearest}
\BIBentryALTinterwordspacing
E.~S. Johnson, S.~Ahmed, S.~S. Dey, and J.-P. Watson, ``A k-nearest neighbor heuristic for real-time dc optimal transmission switching,'' 2020. [Online]. Available: \url{https://arxiv.org/abs/2003.10565}
\BIBentrySTDinterwordspacing

\bibitem{pineda2024learning}
S.~Pineda, J.~M. Morales, and A.~Jim{\'e}nez-Cordero, ``Learning-assisted optimization for transmission switching,'' \emph{Top}, vol.~32, no.~3, pp. 489--516, 2024.

\bibitem{yang2019line}
Z.~Yang and S.~Oren, ``Line selection and algorithm selection for transmission switching by machine learning methods,'' in \emph{2019 IEEE Milan PowerTech}.\hskip 1em plus 0.5em minus 0.4em\relax IEEE, 2019, pp. 1--6.

\bibitem{dey2022node}
S.~S. Dey, B.~Kocuk, and N.~Redder, ``Node-based valid inequalities for the optimal transmission switching problem,'' \emph{Discrete Optimization}, vol.~43, p. 100683, 2022.

\bibitem{hedman2012flexible}
K.~W. Hedman, S.~S. Oren, and R.~P. O’Neill, ``Flexible transmission in the smart grid: optimal transmission switching,'' \emph{Handbook of networks in power systems I}, pp. 523--553, 2012.

\bibitem{gurobi}
\BIBentryALTinterwordspacing
{Gurobi Optimization, LLC}, ``{Gurobi Optimizer Reference Manual},'' 2022. [Online]. Available: \url{https://www.gurobi.com}
\BIBentrySTDinterwordspacing

\bibitem{camm1990cutting}
J.~D. Camm, A.~S. Raturi, and S.~Tsubakitani, ``Cutting big m down to size,'' \emph{Interfaces}, vol.~20, no.~5, pp. 61--66, 1990.

\bibitem{crema2014mathematical}
A.~Crema, ``Mathematical programming approach to tighten a big-m formulation,'' \emph{Escuela de Computacion, Facultad de Ciencias, Universidad Central de Venezuela}, vol.~8, 2014.

\bibitem{klotz2013practical}
E.~Klotz and A.~M. Newman, ``Practical guidelines for solving difficult linear programs,'' \emph{Surveys in Operations Research and Management Science}, vol.~18, no. 1-2, pp. 1--17, 2013.

\bibitem{neumaier2004safe}
A.~Neumaier and O.~Shcherbina, ``Safe bounds in linear and mixed-integer linear programming,'' \emph{Mathematical Programming}, vol.~99, pp. 283--296, 2004.

\bibitem{blumsack2006network}
S.~Blumsack, \emph{Network topologies and transmission investment under electric-industry restructuring}.\hskip 1em plus 0.5em minus 0.4em\relax Pittsburgh, Pennsylvania: Carnegie Mellon University, 2006.

\bibitem{aguilar2025graph}
M.~Aguilar-Moreno, S.~Pineda, and J.~M. Morales, ``A graph-based iterative strategy for solving the all-line transmission switching problem,'' \emph{arXiv preprint arXiv:2502.10333}, 2025.

\bibitem{Oasys2025repository}
{OASYS research group}, ``Sweet spot bound tightening,'' \url{https://github.com/groupoasys/SweetSpotBT}, 2025, accessed: 2025-07-22.

\bibitem{lodi2013performance}
A.~Lodi and A.~Tramontani, ``Performance variability in mixed-integer programming,'' in \emph{Theory driven by influential applications}.\hskip 1em plus 0.5em minus 0.4em\relax INFORMS, 2013, pp. 1--12.

\bibitem{le2015important}
P.~Le~Bodic and G.~L. Nemhauser, ``How important are branching decisions: Fooling mip solvers,'' \emph{Operations Research Letters}, vol.~43, no.~3, pp. 273--278, 2015.

\bibitem{wojtaszek2010faster}
D.~T. Wojtaszek and J.~W. Chinneck, ``Faster mip solutions via new node selection rules,'' \emph{Computers \& operations research}, vol.~37, no.~9, pp. 1544--1556, 2010.

\end{thebibliography}

\end{document}